\documentclass{amsart}
\usepackage{amsmath, amsthm}

\title{ON $(g,\varphi)$-CONTRACTION IN FUZZY METRIC SPACES}

\author{Ioan Golet}
\address{Department of Mathematics "Politechnica" University of Timisoara, 
R-1900 Romania.}
\email{igolet@etv.utt.ro}
\author{Mohd. Rafi}
\address{Faculty of Engineering and Computer Science, The University of Nottingham Malaysia Campus, Jalan Broga, 43500 Semenyih, Selangor Darul Ehsan, Malaysia.}
\email{mdrafzi@yahoo.com}

\subjclass{47H10, 54H25}
\keywords{Generalized fixed point, fuzzy contraction, fuzzy metric space}

\newtheoremstyle{theorem}
  {10pt}		  
  {10pt}  
  {\sl}  
  {\parindent}     
  {\bf}  
  {. }    
  { }    
  {}     
\theoremstyle{theorem}
\newtheorem{theorem}{Theorem}
\newtheoremstyle{lemma}
  {10pt}		  
  {10pt}  
  {\sl}  
  {\parindent}     
  {\bf}  
  {. }    
  { }    
  {}     
\theoremstyle{lemma}
\newtheorem{lemma}[theorem]{Lemma}

\newtheoremstyle{defi}
  {10pt}		  
  {10pt}  
  {\rm}  
  {\parindent}     
  {\bf}  
  {. }    
  { }    
  {}     
\theoremstyle{defi}
\newtheorem{definition}[theorem]{Definition}

\newtheoremstyle{remark}
  {10pt}		  
  {10pt}  
  {\rm}  
  {\parindent}     
  {\bf}  
  {. }    
  { }    
  {}     

\theoremstyle{remark}
\newtheorem{remark}[theorem]{Remark}

\theoremstyle{example}
\newtheorem{example}[theorem]{Example}

\begin{document}

\maketitle

\begin{abstract}
In this paper, we give a generalization of Hicks type contractions and Golet type contractions on fuzzy metric spaces. We prove some fixed point theorems for this new type contraction mappings on fuzzy metric spaces.
\end{abstract}

\section{Introduction and preliminaries}
The notion of fuzzy sets was introduced by Zadeh \cite{Zad}. From then, various concepts
of fuzzy metric spaces were considered in \cite{Geo, Kal, Mih}. Many authors have studied fixed point theory in fuzzy metric spaces. The most interesting references are \cite{Fa1, Fa2, Gre, Had, Pap, Raz}. Some works on intuitionistic fuzzy metric/normed spaces has been carried out intensively in \cite{Raf, Sad}.
In the sequel, we shall adopt usual terminology, notation and convensions of fuzzy
metric spaces introduced by George and Veeramani \cite{Geo}.

\begin{definition}A binary operation $\ast \colon [0,1] \times [0,1]\to [0,1]$ is a continuous $t$-norm if $([0,1],\ast)$ is a topological monoid with unit 1 such that $a \ast b \leq c \ast d$ whenever $a \leq c$ and $b \leq d$ for $a,b,c,d \in [0,1]$.
\end{definition}

\begin{definition}A {\em fuzzy metric space}(briefly $FM$-space) is a triple $(X,M,\ast)$ where $X$ is an arbitrary set, $\ast$ is a continuous $t$-norm and $M\colon X\times X\times [0,+\infty]\to [0,1]$ is a (fuzzy) mapping satisfied the following conditions:\\
For all $x,y,z\in X$ and $s,t>0$,\\
(FM1)  $M(x,y,0)=0$;\\
(FM2)  $M(x,y,t)>0$;\\ 
(FM3)  $M(x,y,t)=1$ if and only if $x=y$;\\
(FM4)  $M(x,y,t)=M(y,x,t)$;\\
(FM5)  $M(x,z,t+s)\geq M(x,y,t)\ast M(y,z,s)$;\\
(FM6)  $M(x,y,\cdot)\colon (0,+\infty)\to [0,1]$ is continuous for all ${x,y} \in X$.
\end{definition}

\begin{lemma}$M(x,y,\cdot)$ is a nondecreasing function for all ${x,y}\in X$.
 \end{lemma}

Let $U_M$ denote the $M$-uniformity, the uniformity generated by fuzzy metric $M$. Then the family

\[ \{U_{\epsilon,\lambda}\colon\epsilon > 0,\lambda \in (0,1)\}\] where
\[U_{\epsilon,\lambda}=\{(x,y)\in X\times X \colon M(x,y,\epsilon)>1-\lambda\},\]
and
\[\{U_{\lambda}\colon \lambda > 0\}\] where \[U_{\lambda}=\{(x,y)\in X \times X \colon M(x,y,\lambda)>1-\lambda\}\]

are bases for this uniformity. 

\begin{remark}Every continuous $t$-norm $\ast$ satisfies \[\sup_{a<1}(a\ast a)=1\] to ensure the existence of the $M-uniformity$ on $X$.
\end{remark}
\medskip

\begin{lemma}\label{lem} If $\ast$ is a continuous $t$-norm, then if $r\in (0,1)$, there is a $s\in (0,1)$ such that $s\ast s\geq r$.
\end{lemma}
\medskip

\begin{definition}In $FM$-space $(X,M,\ast)$, the mapping $f\colon X\to X$ is said to be {\em fuzzy continuous} at $x_0$ if and only if for every $t>0$, there exists $s>0$ such that
\[M(x_{0},y,s)>1-s\Rightarrow M(fx_{0},fy,t)>1-t.\]
\end{definition}
\medskip

The mapping $f\colon X\to X$ is fuzzy continuous if and only if it is fuzzy continuous at every point in $X$.
\medskip
\begin{theorem}Let $(X,M,\ast)$ be a $FM$-space and $U_{M}$ be the $M$-uniformity induced by the fuzzy metric $M$. Then the sequence $\{x_{n}\}in X$ is said to be {\em fuzzy convergence\/} to $x\in X$ (in short, $x_{n}\to x$) if and only if \[\lim_{n\to\infty}M(x_{n},x,t)=1\] for all $t>0$.
\end{theorem}
\medskip

\begin{definition} A sequence $\{x_n\}$ in a $FM$-space $(X,M,\ast)$ is a fuzzy Cauchy sequence if and only if  \[\lim_{n\to +\infty}M(x_{m},x_{n},t)=1\] for every $m,n,>0$ and $t>0$.
\end{definition}
\medskip

\begin{definition} The $FM$-space $(X,M,\ast)$ is said to be {\em fuzzy compact\/} if every sequence $\{x_{n}\}$ in $X$ has a subsequence $\{x_{n_k}\}$ such that $x_{n_k}\to x\in X$.
\end{definition}
\medskip

A fuzzy metric space in which every fuzzy Cauchy sequence is convergent is called a {\em complete fuzzy metric space\/}.
\medskip

\begin{definition} Let $(X,M,\ast)$ be a $FM$-space and $A$ be a nonempty subset of $X$. The fuzzy closure of $A$, denoted by $\overline A$ is the set
\[\overline A=\{y\in X\colon \exists x\in A, M(x,y,\epsilon)>1-\lambda, \epsilon >0, \lambda\in (0,1)\}\] 
\end{definition}

\section{Main Results}

\noindent In our fixed point theorem we consider the FM-space  $(X,M,\ast)$ endowed with the $M$-uniformity.

\begin{definition}\label{def1}Let $\Phi$ be the class of all mappings $\varphi\colon\mathbb R^{+}\to\mathbb R^{+}$ ($\mathbb R^{+} = [0,+\infty)$) with the following properties:

\begin{enumerate}
\item $\varphi$ nondecreasing, 
\item	$\varphi$ is right continuous,
\item	$\lim_{n\to +\infty}\varphi ^{n}(t)=0$ for every $t>0$.
\end{enumerate}
\end{definition}
\medskip

\begin{remark}\label{rem1}(a) It is easy to see that under these conditions, the function $\varphi$ satisfies also $\varphi(t)< t$ for all $t>0$ and therefore $\varphi(0)= 0$.\\
(b)  By property (iii), we mean that for every $\epsilon > 0$ and $\lambda\in (0,1)$ there exists an integer $N(\epsilon,\lambda)$ such that $\varphi^{n}(t) \leq min\{\epsilon,\lambda\}$ whenever $n \geq N(\epsilon,\lambda)$. 
\end{remark}
\medskip

In \cite{Gol}, Golet has introduced $g$-contraction mappings in probabilistic metric spaces. In the following definition, we give the $g$-contraction in the fuzzy setting.
\medskip

\begin{definition}\label{def2}Let $f$ and $g$ are mappings defined on a $FM$-space $(X,M,\ast)$ and let suppose that $g$ is bijective. The mapping $f$ is called a fuzzy {\em $g$-contraction\/} if there exists a $k\in (0,1)$ such that for every ${x,y}\in X$, $t>0$,

\begin{eqnarray}\label{add1}
M(gx,gy,t)>1-t \Rightarrow M(fx,fy,kt)>1-kt.
\end{eqnarray}
\end{definition}

\medskip

By considering a mapping $\varphi \in \Phi$ as given in Definition ~\ref{def1}, we generalize the condition (\ref{add1}) in Definition ~\ref{def2} as follow.
\medskip

\begin{definition}\label{def3}
Let $(X,M,\ast)$ be a $FM$-space and $\varphi \in \Phi$. We say that the mapping $f\colon X \to X$  is a fuzzy {\em $(g,\varphi)$-contraction\/} if there exists a bijective function $g\colon X \to X$ such that for every ${x,y}\in X$ and for every $t>0$,

\begin{eqnarray}\label{add2}
M(gx,gy,t)>1-t \Rightarrow M(fx,fy,\varphi (t))>1- \varphi (t).
\end{eqnarray}
\end{definition}
\medskip

\noindent Note that, if $\varphi (t) = kt$ for $k \in (0,1)$, $t>0$, then the condition (\ref{add2}) is actually the  fuzzy $g$-contraction due to Golet \cite{Gol}. If the function $g$ is identity function, then (\ref{add2}) represents fuzzy ($\varphi - H$)- contraction according to Mihet \cite{Mih}. Hence the $(g,\varphi)$-contraction generalizes the Golet and Mihet's contraction principles respectively in fuzzy metric spaces.
\medskip

The following lemma is reproduced from \cite{Gol} to suits our purposes in fuzzy metric spaces.. 
\medskip

\begin{lemma}\label{lem1}Let $g$ be an injective mapping on $(X,M,\ast)$.\\
(i)   If $M^{g}(x,y,t)=M(gx,gy,t)$, then $(X,M^{g},\ast)$ is a fuzzy metric space.\\
(ii)  If $X_{1}=g(X)$ and $(X,M,\ast)$ is a complete $FM$-space then $(X,M^{g},\ast)$
       is also a complete $FM$-space.\\
(iii) If $(X_{1},M,\ast)$ is fuzzy compact then $(X,M^{g},\ast)$ is also fuzzy compact.
\end{lemma}

\begin{proof}The proof of (i) and (ii) are immediate. To prove (iii), let $\{x_{n}\}$ be a sequence in X. Then, for $u_{n}=gx_{n}$, $\{u_{n}\}$ is a sequence in $X_{1}$ for which we can find a convergent subsequence $\{u_{n_{k}}\}$, say $u_{n_{k}}\to u\in X_{1}$ for $k\to +\infty$. Suppose that the sequence $\{y_{n}\}$ and $y$ in $X$, and set $y_{n}=g^{-1}u_{n_{k}}$,  and $y = g^{-1}u$. Then
\[ M^{g}(y_{n_{k}},y,t) = M(gy_{n_{k}},gy,t) = M(u_{n_{k}},u,t) \to 1\]
as $k \to +\infty$ for every $t>0$. This implies that $(X,M^{g},\ast)$ is fuzzy compact.
\end{proof}
\medskip

\begin{lemma}\label{lem2} If $f$ is a fuzzy $(g,\varphi)$-contraction. Then\\
(i)   $f$ is a fuzzy (uniformly) continuous mapping on $(X,M^{g},\ast)$ with values in $(X,M,\ast)$.
(ii)  $g^{-1}\circ f$ is a continuous mapping on $(X,M^{g},\ast)$ with values into itself.
\end{lemma}

\begin{proof}(i)  Let $\{x_{n}\}$ a sequence in $X$ such that $x_{n}\to x$ in $X$. In $(X,M^{g},\ast)$, this implies that \[ \lim_{n\to +\infty}M^{g}(x_{n},x,t)=\lim_{n\to +\infty}M(gx_{n},gx,t)=1, \forall t>0.\]. By the fuzzy $(g,\varphi)$-contraction (\ref{add2}) it follows that \[ \lim_{n\to +\infty}M(fx_{n},fx,t)\geq \lim_{n\to +\infty}M(fx_{n},fx,\varphi (t))=1, \forall t>0.\]. This implies that $f$ is fuzzy continuous.\\
(ii) Note that, since 
\[ \lim_{n\to +\infty}M(gg^{-1}fx_{n},gg^{-1}fx,t)= \lim_{n\to +\infty}M^{g}(g^{-1}fx_{n},g^{-1}fx,t)=1, \forall t>0,\]
this shows that the mapping $g^{-1} \circ f$ defined on $(X,M^{g},\ast)$ with values in itself is fuzzy continuous.
\end{proof}
\medskip

\begin{theorem}\label{main} Let $f$ and $g$ are two function defined on a complete $FM$-space $(X,M,\ast)$. If $g$ is bijective and $f$ is fuzzy $(g,\varphi)$-contraction, then there exists a unique coincidence point $z\in X$ such that $gz=fz$. 
\end{theorem}

\begin{proof} It is obvious that $M^{g}(x,y,t)>1-t$ whenever $t>1$. Hence, we have $M(gx,gy,t)>1-t$. By condition (\ref{add2}), we have $M(fx,fy,\varphi (t))>1-\varphi (t)$. But, 

\begin{eqnarray*}
M(fx,fy,\varphi (t)) & = & M(gg^{-1}fx,gg^{-1}fy,\varphi(t))\\
& = & M^{g}(g^{-1}fx,g^{-1}fy,\varphi(t))\\
& > & 1-\varphi(t).
\end{eqnarray*}

Now, by letting $h=g^{-1}f$, we have
\[ M^{g}(hx,hy,\varphi (t))>1-\varphi(t).\]
By condition (\ref{add2}), we have

\begin{eqnarray*}
M(fhx,fhy, \varphi^{2}(t)) &=& M(gg^{-1}fhx,gg^{-1}fhy,\varphi^{2}(t))\\
&=& M^{g}(g^{-1}fhx,g^{-1}fhy,\varphi^{2}(t))\\
&=& M^{g}(h^{2}x,h^{2}y,\varphi^{2}(t))\\
&>& 1-\varphi^{2}(t).
\end{eqnarray*} 

By repeating this process, we have \[ M^{g}(h^{n}x,h^{n}y,\varphi^{n}(t))>1-\varphi^{n}(t).\]

Since $\lim_{n\to +\infty}\varphi^{n}(t)=0$, then for every $\epsilon >0$ and $\lambda \in (0,1)$, there exists a positive integer $N(\epsilon,\lambda)$ such that $\varphi^{n}(t)\leq min(\epsilon, \lambda)$, whenever $n\geq N(\epsilon,\lambda)$. Furthermore, since $M$ is nondecreasing  we have
\[ M^{g}(h^{n}x,h^{n}y,\epsilon)\geq M^{g}(h^{n},h^{n},\varphi^{n}(t))>1-\varphi^{n}(t)>1-\lambda.\]
Let $x_{0}$ in $X$ to be fixed and let the sequence $\{x_{n}\}$ in $X$ defined recursively by $x_{n+1}=hx_{n}$, or equivalently by $gx_{n+1}=fx_{n}$. Now, consider $x=x_{p}$ and $y=x_{0}$, then from the above inequality we have
\[ M^{g}(x_{n+p},x_{n},\epsilon) = M^{g}(h^{n}x_{p},h^{n}x_{0},\epsilon)>1-\lambda,\]
for every $n\geq N(\epsilon,\lambda)$ and $p\geq 1$. Therefore $\{x_{n}\}$ is a fuzzy Cauchy sequence in $X$. Since $(X,M,\ast)$ is complete, by Lemma \ref{lem1}, $(X,M^{g},\ast)$ is complete and there is a point $z\in X$ such that $x_{n}\to z$ under $M^{g}$. Since by Lemma \ref{lem2}, $h$ is continuous, we have $z=hz$, i.e., $z=g^{-1}fz$, or equivalently $gz=fz$.  For the uniqueness, assume $gw=fw$ for some $x\in X$. Then, for any $t>0$ and using (\ref{add2}) repeatedly, we can show that after $n$ iterates, we have
\[ M(gz,gw,t)>1-t \Rightarrow M(fz,fw,\varphi^{n}(t))>1-\varphi^{n}(t).\]
Thus, we have $\lim_{n\to +\infty}M(fz,fw,\varphi^{n}(t))=1$ which implies that $fz=fw$.
\end{proof}
\medskip

As an example, if we take in consideration that every metric space $(X,d)$ can be made into a fuzzy metric space $(X,M,\ast)$, in a natural way, by setting $M(x,y,t)= \frac{t}{t+d(x,y)}$, for every ${x,y}\in X$, $t>0$ and $a\ast b=ab$ for every $a,b\in [0,1]$, by Theorem ~\ref{main} one obtains the following fixed point theorem for mappings defined on metric spaces.

\begin{example}Let $f$ and $g$ are two mappings defined on a non-empty set $X$ with values in a complete metric space $(X,d)$, $g$ is bijective and $f$ is a $(g,\varphi)$-contraction, that is, there exists a $\varphi \in \Phi$ such that
\[ d(fx,fy) \leq \varphi(d(gx,gy)),\]
for every ${x,y}\in X$, then there exists a unique element $x^{*}\in X$ such that $fx^{*}=gx^{*}$.
\begin{proof} We suppose that $d(x,y)\in [0,1]$. If this is not true, we define the mapping $d_{1}(x,y)=1-e^{-d(x,y)}$, then the pair $(X,d_{1})$ is a metric space and the uniformities defined by the metrics $d$ and $d_{1}$ are equivalent.\\
Now, let suppose that $f$ is a fuzzy $(g,\varphi)$-contraction for some $\varphi\in \Phi$  on $(X,d)$, $t>0$, $M(gx,gy,t)>1-t$.  Then we have $\frac{t}{t+d(gx,gy)}>1-t$. This implies $d(gx,gy)<t$ and consequently, we have, $\varphi(d(gx,gy))<\varphi(t)$ for all $t>0$. Thus, $d(fx,fy)<\varphi(t)$ which implies that $\frac{\varphi(t)}{\varphi(t)+d(fx,fy)}>1-\varphi(t)$, i.e., $M(fx,fy,\varphi(t))>1-\varphi(t)$. So, the mapping $f$ is a fuzzy $(g,\varphi)$-contraction defined on $(X,M,\ast)$ and the conclusion follows by the Theorem ~\ref{main}.
\end{proof}
\end{example}
\medskip

As a multivalued generalization of the notion of $g$-contraction (Definition \ref{def3}), we shall introduce the notion of a fuzzy $(g,\varphi)$-contraction where $\varphi \in \Phi$ for a multivalued mapping.
\medskip

Let $2^{X}$ be the family of all nonempty subsets of $X$ .

\begin{definition}\label{def4}Let $(X,M,\ast)$ be a $FM$-space, $A$ is a nonempty subset of $X$ and $T\colon A \to 2^{X}$. The mapping $T$ is called a fuzzy $(g,\varphi)$-contraction, where $\varphi \in \Phi$  if there exists a bijective function $g\colon X \to X$ such that for every ${x,y}\in X$ and every $t>0$,
\begin{eqnarray}
M^{g}(x,y,t)>1-t \Rightarrow \\
\nonumber\forall u \in (T\circ g)(x), \exists v\in (T\circ g)(y)\colon  M(u,v,\varphi(t))>1-\varphi(t).
\end{eqnarray}
\end{definition}
\medskip

\begin{definition} Let $(X,M,\ast)$ be a $FM$-space, $A$ is a nonempty subset of $X$ and $T\colon A \to 2^{X}$. We say that $T$ is {\em weakly fuzzy demicompact\/} if for every sequence $\{x_{n}\}$ from $A$ such that $x_{n+1}\in Tx_{n}, n\in\mathbb{N}$ and \[\lim_{n\to+\infty}M(x_{n+1},x_{n},t)=1\]
for every $t>0$, there exists a fuzzy convergent subsequence $\{x_{n_{k}}\}$.
\end{definition}
\medskip

By $cl(X)$ we shall denote the family of all nonempty closed subsets of $X$.
\medskip

\begin{theorem}\label{main2}Let $(X,M,\ast)$ be a complete $FM$-space, $g\colon X \to X$ be a bijective function and $T\colon A \to cl(A)$ where $A\in cl(A)$ a fuzzy $(g,\varphi)$-contraction, where $\ varphi \in \Phi$. If $T$ is weakly fuzzy demicompact then there exists at least one element $x\in A$ such that $x\in Tx$.
\end{theorem}

\begin{proof}Let $x_{0}$, $x_{1}\in (T\circ g)(x_{0})$. Let $t>0$. Since $M^{g}(x_{1},x_{0},t)>0$, it follows that
 \[M^{g}(x_{1},x_{0},t)>1-t.\]
 The mapping $T$ is a fuzzy $(g,\varphi)$-contraction, therefore by Definition\ref{def4} there exists $x_{2}\in (T\circ g)(x_{1})$ such that 
\[M(x_{2},x_{1},\varphi(t))>1-\varphi(t).\] 

Since $M(x_{2},x_{1},\varphi(t))=M^{g}(g^{-1}x_{2},g^{-1}x_{1},\varphi(t))$, we have
\[ M^{g}(g^{-1}x_{2},g^{-1}x_{1},\varphi(t))>1-\varphi(t).\]

Similarly, it follows that there exists $x_{3}\in (T\circ g)(x_{2})$ such that
\[M(x_{3},x_{2},\varphi^{2}(t))>1-\varphi^{2}(t).\] 

By repeating the above process, there exists $x_{n}\in (T\circ g)(x_{n-1})$ $(n\geq 4)$ such that
\[M(x_{n},x_{n-1},\varphi^{n-1}(t))=1-\varphi^{n-1}(t), \forall t>0.\]

By letting $n\to +\infty$, \[\lim_{n\to +\infty}M(x_{n},x_{n-1},\varphi^{n-1}(t))=1, \forall t>0.\] 

Further, by Remark \ref{rem1}(ii), we have \[\lim_{n\to +\infty}M(x_{n},x_{n-1},\epsilon)=1, \forall t>0.\]

Since $T$ is weakly fuzzy demicompact from the above limit, there exists a convergent fuzzy subsequence $\{x_{n}\}$ such that $\lim_{k\to +\infty}x_{n_{k}}=x$ for some $x\in A$. Now, we show that $x\in (T\circ g)(x)$. Since $(T\circ g)(x)=\overline{(T\circ g)(x)}$, we shall prove that $x\in \overline{(T\circ g)(x)}$, i.e. for every $\epsilon >0$ and $\lambda\in (0,1)$, there exists $y\in (T\circ g)(x)$ such that \[M(x,y,\epsilon)>1-\lambda.\]
Note that, since $\ast$ is a continuous $t$-norm, by Lemma \ref{lem}, for $\lambda\in (0,1)$ there is a $\delta\in (0,1)$ such that \[(1-\delta)\ast (1-\delta)\geq 1-\lambda.\]

Further, if $\delta_{1}\in (0,1)$ is such that \[(1-\delta_{1})\ast (1-\delta_{1})\geq 1-\delta\]
and $\delta_{2}= min\{\delta,\delta_{1}\}$, we have
\begin{eqnarray*}
(1-\delta_{2})\ast [(1-\delta_{2})\ast (1-\delta_{2})] &\geq& (1-\delta)\ast [(1-\delta_{1})\ast (1-\delta_{1})]\\
&\geq& (1-\delta)\ast (1-\delta)\\
&>& 1-\delta.
\end{eqnarray*}
Since $\lim_{k\to +\infty}x_{n_{k}}=x$, there exists an integer $k_{1}$ such that \[M(x,x_{n},\epsilon/3)>1-\delta_{2}, \forall k\geq k_{1}.\].
Let $k_{2}$ be an integer such that
\[M(x_{n},x_{n+1},\epsilon/3)>1-\delta_{2}, \forall k\geq k_{2}.\]

Let $s>0$ be such that $\varphi(s)< min\{\epsilon/3,\delta_{2}\}$ and $k_{3}$ be an integer such that \[M^{g}(x_{n_{k}},x,s)>1-s, \forall k\geq k_{3}.\] 

 Since $T$ is $(g,\varphi)$-contraction there exists $y\in (T\circ g)(x)$ such that
\[M(x_{n_{k}+1},y,\varphi(s))>1-varphi(s),\]
and so
\[ M(x_{n_{k}+1},y,\epsilon /3)
\geq (M(x_{n_{k}+1},y,\varphi(s))
>1-\varphi(s)>1-\delta_{2}\]

for every $k\geq k_{3}$. If $k\geq max\{k_{1},k_{2},k_{3}\}$, we have

\begin{eqnarray*}
M(x,y,\epsilon) &\geq& M(x,x_{n_{k}},\epsilon /3)\ast (M(x_{n_{k}}, x_{n_{k}+1},\epsilon /3)\ast M(x_{n_{k}+1},y, \epsilon /3)\\
&\geq& (1-\delta_{2})\ast((1-\delta_{2})\ast (1-\delta_{2}))\\
&>& 1-\lambda.
\end{eqnarray*}
Hence  $x\in\overline{(T\circ g)(x)}=(T\circ g)(x)$. The proof is completed.
\end{proof}
Note that, Theorem \ref{main2} above generalizes the theorem proved by Pap et al in \cite{Pap}.   

\begin{remark}
We note that Theorem \ref{main} can be obtained in a more general setting, namely when $g$ does not satisfy the injective conditions. In this case we can use a pseudo-inverse of $g$ defined as a selector of the multi-valued inverse of $g$. We will consider this situation in our next paper.
\end{remark}

\end{document}